# Balakrishnan Alpha Skew Normal Distribution: Properties and Applications

P. J. Hazarika, S. Shah and S. Chakraborty

Department of Statistics, Dibrugarh University

Dibrugarh, Assam, 786004, India

(Version I: 18th June, 2019)

## Abstract

In this paper a new type of alpha skew distribution is proposed under Balakrishnan Mechanism and some of its related distributions are investigated. The moments and distributional properties are also studied. Suitability of the proposed distribution is tested by conducting data fitting experiments and model adequacy is checked via AIC, BIC in comparison with some related distributions. Likelihood ratio test is carried out to discriminate between normal and proposed distribution.

**Math Classification:** 60E05, 62E10

## 1. Introduction

One cannot undermine the application and the value of normal distribution in real life to model the symmetric data. Now there are many real life situations which seem to be symmetric but due to influences of other factors they depart from symmetry (for details see Chakraborty and Hazarika (2011), and Chakraborty et al., 2015). To tackle this situations Azzalini (1985) introduced the skew-normal distribution, as a natural extension of normal distribution by inducting an additional parameter to introduce asymmetry. A continuous random variable (r.v.) $Z$ follows skew normal (SN) distribution i.e. $Z \sim SN(\lambda)$ if it has probability density function (pdf)

$$f_Z(z;\lambda) = 2\phi(z)\Phi(\lambda z); \quad z, \lambda \in R \qquad (1)$$

where, $\phi$ and $\Phi$ are the pdf and cumulative distribution function (cdf) of standard normal distribution respectively.

*Preprint*

Balakrishnan (2002) as a discussant in Arnold and Beaver (2002) proposed the generalization of the skew normal density and studied its properties. The pdf of the same distribution is

$$f_Z(z; \lambda, n) = \phi(z)[\Phi(\lambda z)]^n / C_n(\lambda); \quad z, \lambda \in R \qquad (2)$$

where, $n$ is positive integer and $C_n(\lambda) = E(\Phi^n(\lambda U))$, $U \sim N(0,1)$. In particular, if $\lambda = 1$ the Balakrishnan skew normal density becomes skew normal density of Azzalini (1985). Furthermore, Sharafi and Behbodian (2008) extensively studied its different forms and properties. Bahrami et al. (2009) introduced the two parameter Balakrishnan skew normal distribution. Yadegari et al. (2008) discussed the generalization of Balakrishnan skew normal distribution

Huang and Chen (2007) proposed the general formula for the construction of skew-symmetric distributions starting from a symmetric (about 0) pdf $h(.)$ by introducing the concept of skew function $G(.)$, a Lebesgue measurable function such that, $0 \leq G(z) \leq 1$ and $G(z) + G(-z) = 1$, $z \in R$, almost everywhere. A random variable $Z$ is said to be skew symmetric if the pdf is of the following form

$$f(z) = 2h(z)G(z); z \in R \qquad (3)$$

Olivero in 2010 developed a new form of skew distribution which exhibits both unimodal as well as bimodal behavior and named it as alpha skew normal distribution. An r.v. $Z$ is said to follow alpha skew normal distribution ($\text{ASN}(\alpha)$) if its pdf is given by

$$f(z; \alpha) = \{(1-\alpha z)^2 + 1\}\varphi(z)/(2+\alpha^2); z, \alpha \in R \qquad (4)$$

Using the same approach of Olivero (2010), Harandi and Alamatsaz (2013) and Hazarika and Chakraborty (2014) respectively, investigated a class of alpha skew Laplace distributions and alpha skew logistic distributions. Venegas et al. (2016) and Louzada et al. (2017) studied the logarithmic form and bivariate form of alpha-skew-normal distribution, respectively. Sharafi et al. (2017) discussed the generalization of alpha-skew-normal distribution.

In this article a new version of alpha skew normal distribution is proposed by considering methodology advocated by Balakrishnan in 2002 and some of its basic properties are investigated.

## 2. Balakrishnan Alpha Skew Normal Distribution

In this section we introduce the generalized version of bimodal skew normal distribution of Olivero (2010) and proposed Balakrishnan alpha skew normal distribution.

**Definition 1:** A random variable $Z$ is said to follow generalized bimodal normal distribution if its pdf is given by

$$f(z) = \frac{z^n}{C}\varphi(z); \quad z \in R \tag{5}$$

where, $n$ is positive even integers and $C$ is normalizing constant. Symbolically, we can write $Z \sim BN(n)$. The shapes of pdfs with different choices of $n$ are shown in figure 1.

**Remark 1:** The pdf in (5) has at most two modes and the same has been seen from the equation $f'(z) = \frac{x^{n-1}(x^2 - n)}{C}\varphi(z) = 0$. This equation has only three zero, therefore the pdf in (5) have only two modes.

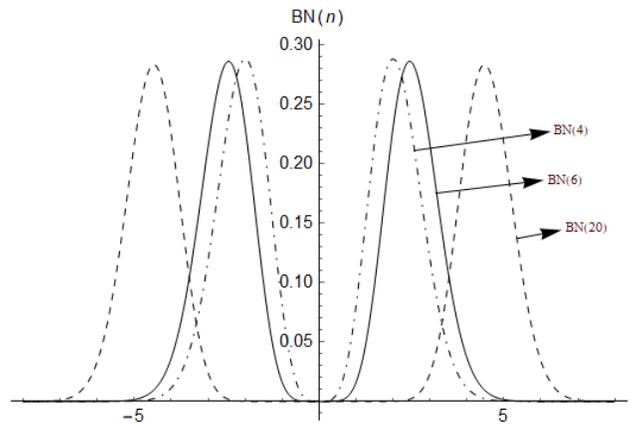

**Figure 1:** Plots of pdf of $BN(n)$

**Definition 2:** A random variable $Z$ with density function

$$f_Z(z;\alpha) = \frac{1}{C_2(\alpha)}\left(\frac{(1-\alpha z)^2 + 1}{2 + \alpha^2}\right)^2 \varphi(z); \quad z, \alpha \in R \tag{6}$$

where, $C_2(\alpha) = 3 - \frac{4}{2 + \alpha^2}$, is said to follow **Balakrishnan alpha skew normal** distribution with parameter $\alpha$.

In the rest of this article, we shall refer the distribution in (6) as $BASN_2(\alpha)$. The plots of pdf are depicted in figure 2 for different choices of the parameter $\alpha$.

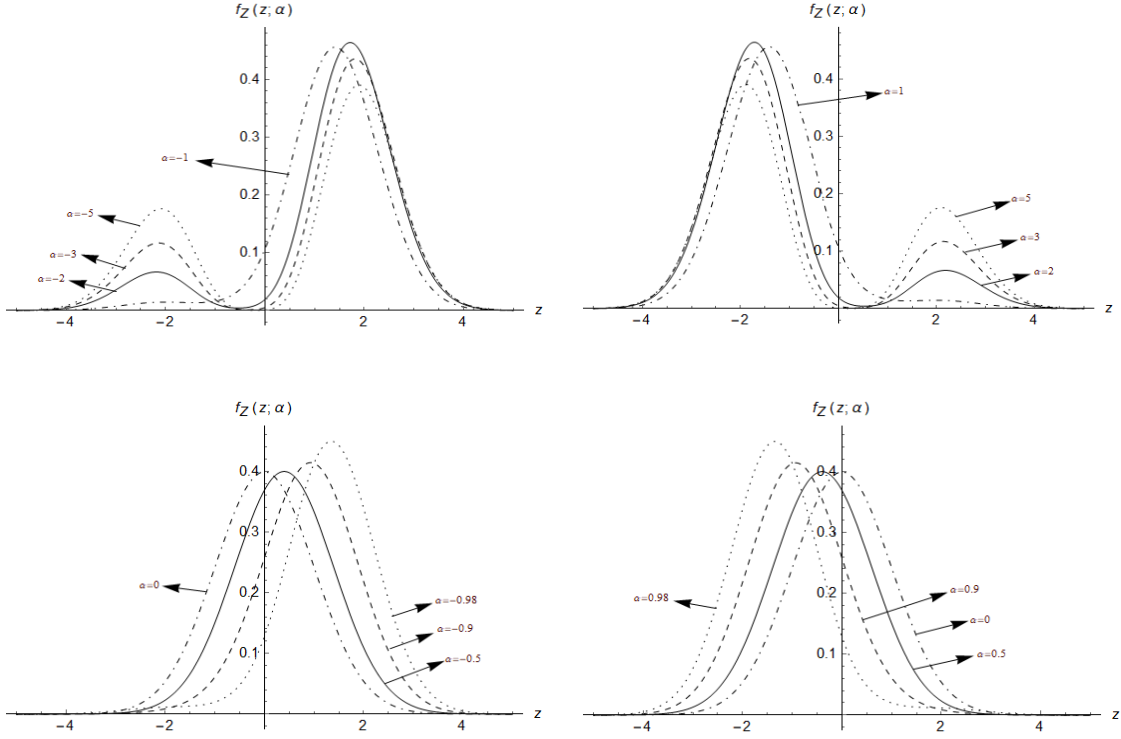

**Figure 2:** Plots of pdf of $BASN_2(\alpha)$

It's obvious to note and check that $BASN_2(\alpha)$ is bimodal when $|\alpha| \geq 1$.

**Remark 2:** The pdf of the proposed alpha skew normal distribution is constructed using the formula (2), by taking $\Phi(.) = \dfrac{(1-\alpha z)^2 + 1}{2+\alpha^2}$ and $n = 2$.

**Properties of** $BASN_2(\alpha)$**:**

i)   $BASN_2(0) = \varphi(z)$

ii)  If $\alpha \to \pm\infty$ then pdf of Z becomes $f_Z(z) = \dfrac{z^4}{3\sqrt{2\pi}} \exp\left(-\dfrac{z^2}{2}\right)$ i.e., $Z \sim BN(4)$

iii) If $Z \sim BASN_2(\alpha)$ then $-Z \sim BASN_2(-\alpha)$

iv)  $BASN_2(\alpha)$ have at most two modes.

*Proof:* To show $BASN_2(\alpha)$ have at most two modes, which is equivalent to prove that the following equation have three zeros.

$$Df_Z(z;\alpha) = \frac{[(1-\alpha z)^2 + 1](\alpha^2 z^3 - 4\alpha^2 z - 2\alpha z^2 + 4\alpha + 2z)\varphi(z)}{C_2(\alpha)(2+\alpha^2)^2} = 0 \qquad (7)$$

It is easy to show that the eqn. (7) has at most three real zeros because $(1-\alpha z)^2 + 1 = 0$ will have two complex roots, $\alpha^2 z^3 - 4\alpha^2 z - 2\alpha z^2 + 4\alpha + 2z = 0$ has three real roots and $\varphi(z) \neq 0$. The same can be depicted from the contour plot of the eqn. (7) given in figure 3. It is also observed that approximately for $-0.95 < \alpha < 0.95$, $BASN_2(\alpha)$ remains unimodal (see the second plot in figure 3).

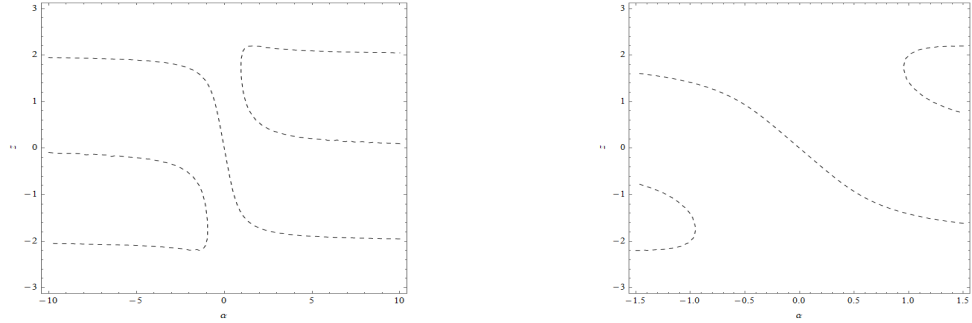

**Figure 3:** Contour plots of the equation $Df_Z(z;\alpha) = 0$

**Proposition 1:** If $Z \sim BASN_2(\alpha)$ then its cdf is given by

$$F_Z(z;\alpha) = \Phi(z) + \frac{\alpha\{8 - 8\alpha z + 4\alpha^2(2+z^2) - \alpha^3 z(3+z^2)\}}{C_2(\alpha)(2+\alpha^2)^2}\varphi(z) \tag{8}$$

where, $\Phi(z)$ is the cdf of standard normal distribution.

*Proof:* see Appendix A.

The plot of cdf with different choices parameter is shown in figure 4.

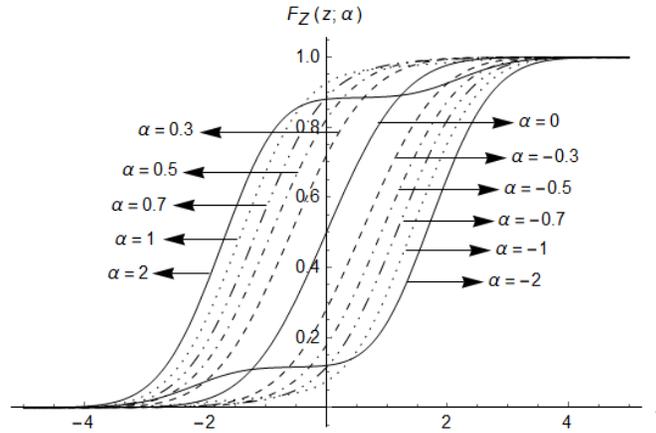

**Figure 4:** Plots of cdf of $BASN_2(\alpha)$

For $-1 < \alpha < 0 \, (0 < \alpha < 1)$, we can say that the standard normal is stochastically smaller (larger) than $BASN_2(\alpha)$ as seen in the figure 4.

**Remark 3:** In particular, if $\alpha \to \pm\infty$ then the cdf of $BASN_2(\alpha)$ becomes the cdf of $BN(4)$ and which is nothing but $F_Z(z) = \Phi(z) - \frac{z(3+z^2)}{3}\varphi(z)$

**Proposition 2:** If $Z \sim BASN_2(\alpha)$ then

$$E(Z^n) = \begin{cases} \dfrac{2^{-\frac{(n+4)}{2}}\left\{\dfrac{\alpha^4(n+4)!}{((n+4)/2)!} + \dfrac{16\alpha^2(n+2)!}{((n+2)/2)!} + \dfrac{16n!}{(n/2)!}\right\}}{C_2(\alpha)(2+\alpha^2)^2}; & \text{when } n \text{ is even} \\[2ex] \dfrac{-2^{\frac{1-n}{2}}\alpha\left\{\dfrac{\alpha^2(n+3)!}{((n+3)/2)!} + \dfrac{4(n+1)!}{((n+1)/2)!}\right\}}{C_2(\alpha)(2+\alpha^2)^2}; & \text{when } n \text{ is odd} \end{cases} \quad (9)$$

*Proof:* see Appendix B.

**Remark 4:** The expression (9) can be rewritten with the help of Gamma function as

$$E(Z^n) = \begin{cases} \dfrac{2^{\frac{n}{2}}\{4+(1+n)\alpha^2(8+(3+n)\alpha^2)\}\Gamma(\frac{1+n}{2})}{C_2(\alpha)\,\Gamma(1/2)(2+\alpha^2)^2}, & \text{when } n \text{ is even} \\[2ex] -\dfrac{2^{\frac{5+n}{2}}\alpha\{2+(2+n)\alpha^2\}\Gamma(1+\frac{n}{2})}{C_2(\alpha)\Gamma(1/2)(2+\alpha^2)^2}, & \text{when } n \text{ is odd} \end{cases}$$

In particular, $E(Z) = \dfrac{-4\alpha}{(2+\alpha^2)}$, $E(Z^2) = 5 - \dfrac{4}{(2+\alpha^2)} - \dfrac{4}{(2+3\alpha^2)}$, $Var(Z) = \dfrac{(2+5\alpha^2)(4+3\alpha^4)}{(2+\alpha^2)^2(2+3\alpha^2)}$

$E(Z^3) = \dfrac{-12\alpha(2+5\alpha^2)}{4+8\alpha^2+3\alpha^4}$ and $E(Z^4) = 35 - \dfrac{48}{(2+\alpha^2)} - \dfrac{16}{(2+3\alpha^2)}$.

**Remark 5:** By optimizing $E(Z)$ and $Var(Z)$ with respect to $\alpha$ we get the following bounds.

i. $-1.414 \leq E(Z) \leq 1.414$

ii. $0.972 \leq Var(Z) \leq 4.7966$

The same can be easily visualized from figure 5 and figure 6.

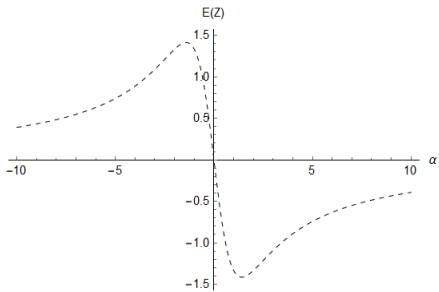
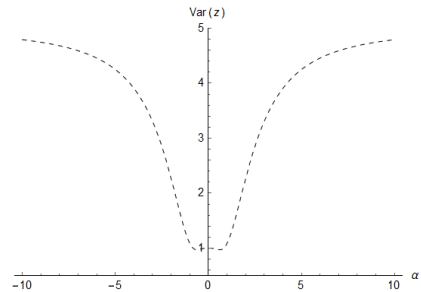

**Figure 5:** Plot of mean          **Figure 6:** Plot of variance

**Remark 6:** The expression for skewness $(\beta_1)$ and kurtosis $(\beta_2)$ are respectively given by

$$\beta_1 = \frac{64\alpha^6(2+3\alpha^2)(4+15\alpha^4)^2}{(8+20\alpha^2+6\alpha^4+15\alpha^6)^3} \text{ and}$$

$$\beta_2 = \frac{3(2+3\alpha^2)(32+112\alpha^2+144\alpha^4+216\alpha^6+410\alpha^8+35\alpha^{10})}{(8+20\alpha^2+6\alpha^4+15\alpha^6)^2}$$

**Remark 7:** By optimizing $\beta_1$ and $\beta_2$ with respect to $\alpha$ we get the following bounds.

  i.  $2.5359 \leq \beta_1 \leq 0$
  ii. $6.7684 \leq \beta_2 \leq 3$

The same can be easily visualized from figure 7 and figure 8.

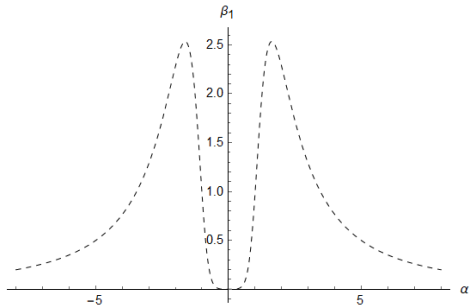
**Figure 7:** Plot of skewness

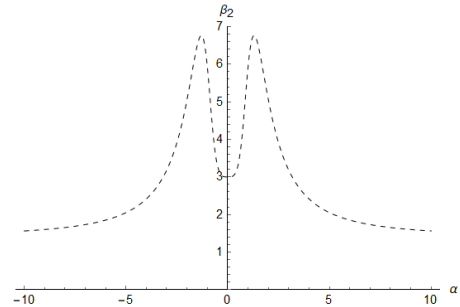
**Figure 8:** Plot of kurtosis

**Proposition 3:** If $Z \sim BASN_2(\alpha)$ then its mgf is given by

$$M_Z(t) == \frac{M_X(t)[\alpha^4 t^4 + 6\alpha^4 t^2 + 3\alpha^4 - 4\alpha^3 t^3 - 34\alpha^3 t + 8\alpha^2 t^2 + 8\alpha^2 - 8\alpha t + 4]}{C_2(\alpha)(2+\alpha^2)^2} \quad (10)$$

where, $M_X(t)$ is the mgf of standard normal variable.

*Proof:* see Appendix C.

**Remark 8:** The $BASN_2(\alpha)$ distribution can be expressed into two component of symmetric and asymmetric part as shown below

$$f(z;\alpha) == \frac{\alpha^4 z^4 + 8\alpha^2 z^2 + 4}{C_2(\alpha)(2+\alpha^2)^2}\varphi(z) + \frac{(-4\alpha^3 z^3 - 8\alpha z)}{C_2(\alpha)(2+\alpha^2)^2}\varphi(z) \quad (11)$$

In equation (11) the 1st part is symmetric and the second part is asymmetric one and the symmetric part is symbolically denoted by $SCBASN_2(\alpha)$. In particular, if $\alpha = 0$ then $SCBASN_2(\alpha)$ becomes standard normal distribution.

**Proposition 4:** If $Z \sim SCBASN_2(\alpha)$ then its cdf is given by

$$F(z) = \Phi(z) - \frac{\alpha^2[\alpha^2 z^3 + 3z\alpha^2 + 8z]}{C_2(\alpha)(2+\alpha^2)^2}\varphi(z) \qquad (12)$$

where, $\Phi(z)$ is the cdf of standard normal distribution.

*Proof:* see Appendix D.

**Proposition 5:** $Z \sim SCBASN_2(\alpha)$ then its mgf is given by

$$M_Z(t) = \frac{[\alpha^4 t^4 + 6\alpha^4 t^2 + 3\alpha^4 + 8\alpha^2 t^2 + 8\alpha^2 + 4]}{C_2(\alpha)(2+\alpha^2)^2} M_X(t) \qquad (13)$$

*Proof:* see Appendix E.

**Remark 9:** To generate the random number $Z$ from $BASN_2(\alpha)$ distribution for different choices of the parameter $\alpha$ one can adopt the acceptance sampling method with the following steps:

    I: Generate random number $U$ from Uniform $(0,1)$

    II: Generate random number $H$ from $SCBASN_2(\alpha)$.

    III: Set $Z = H$ if $U < \frac{1}{\Delta}\frac{f(H)}{f_1(H)}$, otherwise, go to step I and continue the process.

Where, $\Delta = Sup\left[\frac{f(Z)}{f_1(Z)}\right] = \frac{1}{3}(1+2\sqrt{3})$ and $f(.)$ and $f_1(.)$ are pdf of $BASN_2(\alpha)$ and $SCBASN_2(\alpha)$ respectively.

## 3. Half $BASN_2(\alpha)$ Distribution

A half Balakrishnan-alpha-skew normal $HBASN_2(\alpha)$ distribution truncated below '0' is given by

$$f_T(t;\alpha) = \frac{[(1-\alpha t)^2 + 1]^2}{3\alpha^4 - 8\alpha^3 b + 8\alpha^2 - 8\alpha b + 4}\psi(t); \quad t > 0 \qquad (14)$$

Where, $\psi(t)$ is the pdf of the standard half-normal distribution and $b = \sqrt{\frac{2}{\pi}}$.

This can be considered as a potential life time distribution. The corresponding survival function $S_T(t;\alpha)$ and the hazard rate functions $h_T(t;\alpha)$ of $TBASN_2(\alpha)$ can be expressed as below

$$S_T(t;\alpha) = \frac{\psi(t)\alpha(8-8\alpha t+8\alpha^2+4\alpha^2 t^2-3\alpha^3 t-\alpha^3 t^3)-C_2(\alpha)(2+\alpha^2)^2\overline{\Psi}(t)}{-4+\alpha[8c+\alpha(-8+8\alpha b-3\alpha^2)]}$$

and $h_T(t;\alpha) = \dfrac{[(1-\alpha t)^2+1]^2}{\alpha(-8+8\alpha t-8\alpha^2-4\alpha^2 t^2+3\alpha^3 t+\alpha^3 t^3)+C_2(\alpha)(2+\alpha^2)^2\psi(t)\overline{\Psi}(t)}$.

where, $\Psi(t)$ and $\overline{\Psi}(t)$ are respectively the cdf and survival function of the standard half-normal distribution.

We have plotted the $h_Z(t;\alpha)$ for the suitable values of the parameter $\alpha$, in figure 9 to study its behavior graphically.

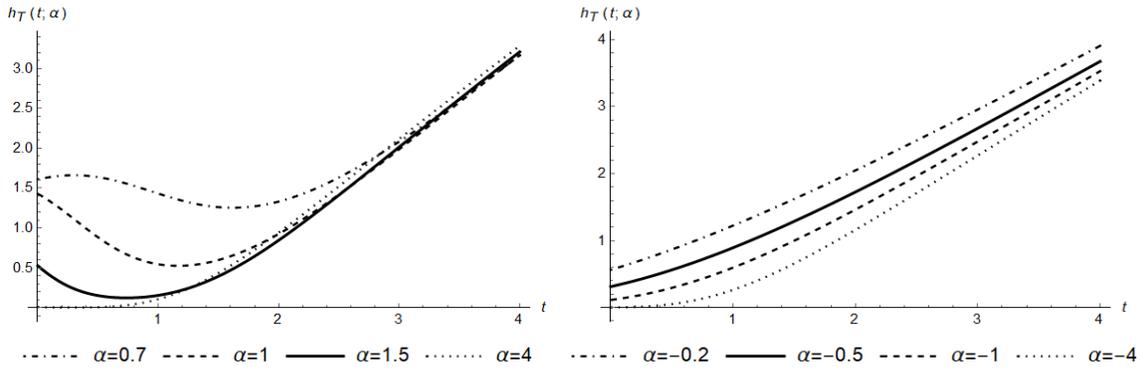

**Figure 9:** Plots of hazard rate function of $HBASN_2(\alpha)$

It can be observed from figure 9 that the hazard rate is increasing for $\alpha \leq 0$ while it assumes bathtub shape for $\alpha > 0.7$. For the values of $0 < \alpha \leq 0.7$, the hazard rate takes different shapes. Therefore, the hazard rate function of $HBASN_2(\alpha)$ distribution assumes different useful shapes depending on the choice of the values of the parameter $\alpha$, and thus has the potential to be a flexible life time model.

**Remark 10:** In particular for $\alpha = 0$, $HBASN_2(\alpha)$ distribution reduces to standard half-normal distribution.

## 4. Some Extensions of $BASN_2(\alpha)$ Distribution

### 4.1. The Bivariate $BASN_2(\alpha)$ Distribution

**Definition 3:** A random vector $Z = (Z_1, Z_2)$ is two dimensional (bivariate) Balakrishnan alpha skew normal distributed if it has the following density

$$f(z;\alpha_1,\alpha_2,\rho) = \frac{[(1-\alpha_1 z_1-\alpha_2 z_2)^2+1]^2}{C(\alpha_1,\alpha_2,\rho)}\varphi_2(z); z \in R^2,\ \alpha_1,\alpha_2 \in R \qquad (15)$$

Where, $C(\alpha_1, \alpha_2, \rho) = (2 + \alpha_1^2 + 2\rho\alpha_1\alpha_2 + \alpha_2^2)(2 + 3\alpha_1^2 + 6\rho\alpha_1\alpha_2 + 3\alpha_2^2)$, and $\varphi_2(z)$ is the pdf of a bivariate normal distribution $N_2\left(\begin{bmatrix} 0 \\ 0 \end{bmatrix}, \begin{bmatrix} 1 & \rho \\ \rho & 1 \end{bmatrix}\right)$. We denote it by $Z \sim BBASN_2(\alpha_1, \alpha_2, \rho)$

**Special cases of $BBASN_2(\alpha_1, \alpha_2, \rho)$:**

- If $\alpha_1 = \alpha_2 = 0$, then $Z \sim N_2\left(\begin{bmatrix} 0 \\ 0 \end{bmatrix}, \begin{bmatrix} 1 & \rho \\ \rho & 1 \end{bmatrix}\right) = \varphi_2(z)$

- If $\alpha_2 = 0$, then pdf of $Z_1$ is $\dfrac{[(1-\alpha_1 z_1)^2 + 1]^2}{4 + 8\alpha_1^2 + 3\alpha_1^4}\varphi_2(z_1)$ and if $\alpha_1 = 0$, then pdf of $Z_2$ is

  $\dfrac{[(1-\alpha_2 z_2)^2 + 1]^2}{4 + 8\alpha_2^2 + 3\alpha_2^4}\varphi_2(z_2)$

- If $\alpha_1 = \alpha_2 \to \pm\infty$, then $f(z; \alpha_1, \alpha_2, \rho) = \dfrac{(z_1 + z_2)^4}{12(1+\rho)^2}\varphi_2(z); \quad z \in R$

- If $Z \sim BBASN_2(\alpha_1, \alpha_2, \rho)$, then $-Z \sim BBASN_2(-\alpha_1, -\alpha_2, \rho)$

## 4.2 A Two-parameter $BASN_2(\alpha)$ Distribution

**Definition 4:** A random variable $Z$ has a two-parameter Balakrishnan alpha skew normal distribution with parameters $\alpha_1, \alpha_2 \in R$, denoted by $TPBASN_2(\alpha_1, \alpha_2)$, if its pdf is

$$f(z; \alpha_1, \alpha_2) = \dfrac{\Phi^2(\alpha_1 z)\Phi^2(\alpha_2 z)}{C(\alpha_1, \alpha_2)}\varphi(z); \quad z \in R \tag{16}$$

where, $C(\alpha_1, \alpha_2) = 32\alpha_1\alpha_2(2 + 3\alpha_2^2) + 48\alpha_1^3\alpha_2(2 + 5\alpha_2^2) + 4(4 + 8\alpha_2^2 + 3\alpha_2^4) + 8\alpha_1^2(4 + 24\alpha_2^2 + 15\alpha_2^4) + 3\alpha_1^4(4 + 40\alpha_2^2 + 35\alpha_2^4)$

$\Phi^2(\alpha_1 z) = [(1-\alpha_1 z)^2 + 1]^2$; and $\Phi^2(\alpha_2 z) = [(1-\alpha_2 z)^2 + 1]^2$

**Special cases of $TPBASN_2(\alpha_1, \alpha_2)$:**

- If $\alpha_1 = \alpha_2 = 0$, then $Z \sim N(0,1) = \varphi(z)$

- If $\alpha_2 = 0$, then $Z \sim BASN_2(\alpha_1)$ and if $\alpha_1 = 0$, then $Z \sim BASN_2(\alpha_2)$

- If $\alpha_1 = \alpha_2 = \alpha$, then $f(z; \alpha) = \dfrac{((1-\alpha z)^2 + 1)^4}{16 + 128\alpha^2 + 408\alpha^4 + 480\alpha^6 + 105\alpha^8}\varphi(z) = BASN_4(\alpha)$

- If $\alpha_1 = \alpha_2 \to \pm\infty$, then $f(z) = \dfrac{z^8}{105}\varphi(z)$

- If $Z \sim TPBASN_2(\alpha_1, \alpha_2)$, then $-Z \sim TPBASN_2(-\alpha_1, -\alpha_2)$

## 4.3 Balakrishnan Alpha-Beta Skew Normal Distribution

**Definition 5:** If the density of the random variable $Z$ has pdf given by

$$f(z;\alpha,\beta) = \frac{[(1-\alpha z - \beta z^3)^2 + 1]^2}{C(\alpha,\beta)} \varphi(z); \quad z,\alpha,\beta \in R \tag{17}$$

Then, we say that $Z$ is distributed according to the Balakrishnan alpha-beta skew normal distribution with parameter $\alpha$ and $\beta$.

where, $C(\alpha,\beta) = 4 + 3\alpha^4 + 60\alpha^3\beta + 12\alpha\beta(4+315\beta^2) + \alpha^2(8+630\beta^2) + 15\beta^2(8+693\beta^2)$. We denote it by $Z \sim BABSN_2(\alpha)$.

**Special cases of $BABSN_2(\alpha)$:**

- If $\beta = 0$, then we get $BASN_2(\alpha)$ distribution and is given by

$$f(z) = \frac{[(1-\alpha z)^2 + 1]^2}{C_2(\alpha)(2+\alpha^2)^2} \varphi(z)$$

- If $\alpha = 0$, then we get

$$f(z) = \frac{[(1-\beta z^3)^2 + 1]^2}{(4+15\beta^2(8+693\beta^2))} \varphi(z)$$

This equation is known as Balakrishnan beta skew normal ($BBSN_2(\alpha)$) distribution.

- If $\alpha = \beta = 0$, then we get the standard normal distribution

- If $\alpha \to \pm\infty$, then we get the bimodal normal ($BN(4)$) distribution given by

$$f(z) = \frac{z^4}{3} \varphi(z)$$

- If $\beta \to \pm\infty$, then we get the bimodal normal ($BN(12)$) distribution given by

$$f(z) = \frac{z^{12}}{10395} \varphi(z)$$

- If $Z \sim BABSN_2(\alpha,\beta)$, then $-Z \sim LBABSN_2(-\alpha,-\beta)$.

## 4.4 Generalization of $BASN_2(\alpha)$ Distribution

**Definition 6:** If the density of the random variable $Z$ has pdf given by

$$f(z;\alpha,\lambda) = \frac{[(1-\alpha z)^2 + 1]^2}{C(\alpha,\lambda)} \varphi(z)\Phi(\lambda z); \quad z,\alpha \in R, \lambda > 0 \tag{18}$$

then, we say that $Z$ is distributed according to the Generalized Balakrishnan alpha skew normal distribution with parameter $\alpha$ and $\lambda$.

where, $C(\alpha,\lambda) = (2 + 4\alpha^2 + 1.5\alpha^4) - b(2\alpha^3(1-\delta^2)\delta + 4\alpha\delta + 4\alpha^3\delta)$; $\delta = \dfrac{\lambda}{\sqrt{1+\lambda^2}}$; $b = \sqrt{\dfrac{2}{\pi}}$;

$\varphi(z)$ and $\Phi(\lambda z)$ are the pdf and cdf of the standard normal distribution. We denote it by $Z \sim GBASN_2(\alpha,\lambda)$.

**Special cases of $GBASN_2(\alpha,\lambda)$:**

- If $\alpha = 0$, then $Z \sim SN(\lambda)$
- If $\lambda = 0$, then $Z \sim BASN_2(\alpha)$
- If $\alpha = \lambda = 0$, then $Z \sim N(0,1)$
- If $\alpha \to \pm\infty$, then $f(z;\alpha,\lambda) \to \dfrac{2z^4}{3}\varphi(z)\Phi(\lambda z)$
- If $\lambda \to +\infty$, then $f(z;\alpha,\lambda) \to \dfrac{[(1-\alpha z)^2 + 1]^2}{(2 + 4\alpha^2 + 1.5\alpha^4) - b(2\alpha^3 + 4\alpha + 4\alpha^3)}\varphi(z)I(x>0)$,

  and if $\lambda \to -\infty$, then $f(z;\alpha,\lambda) \to \dfrac{[(1-\alpha z)^2 + 1]^2}{(2 + 4\alpha^2 + 1.5\alpha^4) - b(2\alpha^3 + 4\alpha + 4\alpha^3)}\varphi(z)I(x<0)$,

  where, $b = \sqrt{\dfrac{2}{\pi}}$ and $I(.)$ is an indicator function.

- If $Z \sim GBASN_2(\alpha,\lambda)$, then $-Z \sim GBASN_2(-\alpha,-\lambda)$

## 4.5 The Log-$BABSN_2(\alpha)$ Distribution

In this section, using the work of (Venegas et al., 2016), we present the definition and some simple properties of log-Balakrishnan alpha skew normal distribution.

Let $Z = e^Y$, then $Y = Log(Z)$, therefore, the density function of Z is defined as follows:

**Definition 7:** If the density of the random variable $Z$ has pdf given by

$$f(z;\alpha) = \dfrac{[(1-\alpha y)^2 + 1]^2}{C_2(\alpha)(2+\alpha^2)^2}\varphi(y); \quad z>0, \alpha \in R \qquad (19)$$

then, we say that Z is distributed according to the log-Balakrishnan alpha skew normal distribution with parameter $\alpha$. Where, $y = Log(z)$ and $\varphi(y)$ is the pdf of the standard log-normal distribution. We denote it by $Z \sim LBASN_2(\alpha)$.

**Special cases of $LBASN_2(\alpha)$:**

- If $\alpha = 0$, then we get the standard log-normal distribution given by

$$f(z) = \dfrac{\varphi(y)}{z}$$

- If $\alpha \to \pm\infty$, then we get the log-bimodal normal $LBN(4)$ distribution given by

$$f_Z(z) = \frac{y^4}{3z}\varphi(y)$$

- If $Z \sim LBASN_2(\alpha)$, then $-Z \sim LBASN_2(-\alpha)$.

## 5. Parameter Estimation of $BASN_2(\alpha)$

If $Z \sim BASN_2(\alpha)$ distribution then $Y = \mu + \sigma Z$ is the location ($\mu$) and scale ($\sigma$) extension of Z and has the pdf is given by

$$f_Y(y;\mu,\sigma,\alpha) = \frac{1}{C_2(\alpha)}\left(\frac{\left\{1-\alpha\left(\frac{y-\mu}{\sigma}\right)\right\}^2+1}{2+\alpha^2}\right)^2 \frac{e^{-\frac{1}{2}\left(\frac{y-\mu}{\sigma}\right)^2}}{\sigma\sqrt{2\pi}} \; ; \; y,\mu,\alpha \in \mathbb{R} \text{ and } \sigma > 0 \quad (20)$$

Symbolically, we write as $Y \sim BASN_2(\alpha,\mu,\sigma)$.

### 5.1. Method of Moments

Let $y_1, y_2, ..., y_n$ be a random sample of size $n$ drawn from $BASN_2(\alpha,\mu,\sigma)$ distribution in eqn.(21) and $m_1, m_2$ and $m_3$ are respectively the first three sample raw moments. Then, the moment estimates of the three parameters $\mu, \sigma$ and $\alpha$ are obtained by

$$m_1 = \mu - \frac{4\alpha\sigma}{2+\alpha^2} \Rightarrow \mu = \frac{4\alpha\sigma}{2+\alpha^2} + m_1 \quad (21)$$

$$m_2 = \mu^2 - \frac{\sigma(16\alpha\mu + 24\alpha^3\mu) - \sigma^2(4\sigma^2 + 24\alpha^2\sigma^2 + 15\alpha^4)}{(2+\alpha^2)(2+3\alpha^2)} \quad (22)$$

$$m_3 = \mu^3 + \frac{-3\sigma(8\alpha\mu^2 + 12\alpha^3\mu^2 - 4\mu\sigma - 24\alpha^2\mu\sigma - 15\alpha^4\mu\sigma + 8\alpha\sigma^2 + 20\alpha^3\sigma^2)}{4+8\alpha^2+3\alpha^4} \quad (23)$$

Substituting the value of $\mu$ from eqn. (21) in eqn. (22) and solving for $\sigma^2$, we get

$$\sigma^2 = \frac{(m_2 - m_1^2)(2+\alpha^2)^2(2+3\alpha^2)}{(8+20\alpha^2+6\alpha^4+15\alpha^6)} \quad (24)$$

Finally, by putting these values of $\mu$ and $\sigma^2$ in eqn.(23), we get the following equation in $\alpha$

$$m_3 = \frac{\begin{bmatrix}[C_2(\alpha)d_1^2(d_3d_4m_1 - 4\alpha d_1d_2(m_1^2 - m_2))^3 + 12\alpha\, d_1^2\, d_2^2(d_3d_4m_1 - 4\alpha d_1d_2(m_1^2 - m_2))^2(m_1^2 - m_2) + \\ 3d_1^4\, d_2^2 d_5(d_3d_4m_1 - 4\alpha d_1d_2(m_1^2 - m_2))(m_1^2 - m_2)^2 + 12\alpha d_1^6\, d_2^3 d_3(m_1^2 - m_2)^3]\end{bmatrix}}{C_2(\alpha)d_1^2 d_3^3 d_4^3}$$

(25)

where, $d_1 = 2+\alpha^2$; $d_2 = 2+3\alpha^2$; $d_3 = 2+5\alpha^2$; $d_4 = 4+3\alpha^4$; and $d_5 = 4+24\alpha^2+15\alpha^4$

Furthermore, the value of $\alpha$ is estimated numerically as the exact solution of the eqn.(25) is not easily tractable. Once $\alpha$ is estimated, the rest of the two parameters namely, $\mu$ and $\sigma$ can be estimated directly from eqn. (21) and eqn. (24) respectively.

## 5.2. Maximum likelihood Method

*Likelihood function:*

Let $y_1, y_2, ..., y_n$ be a random sample of size $n$ drawn from $BASN_2(\alpha, \mu, \sigma)$ distribution of eqn. (20), then the log-likelihood function for $\theta = (\alpha, \mu, \sigma)$ is given by

$$l(\theta) = 2\sum_{i=1}^{n}\log\left[\left\{1-\alpha\left(\frac{y_i-\mu}{\sigma}\right)\right\}^2+1\right] - n\log(2+\alpha^2) - n\log\sigma - n\log(2+3\alpha^2) - \frac{n}{2}\log(2\pi) - \frac{1}{2}\sum_{i=1}^{n}\left(\frac{y_i-\mu}{\sigma}\right)^2 \qquad (26)$$

Differentiating this eqn. (26) above partially with respect to the parameters $\alpha, \mu,$ and $\sigma$, the following likelihood equations are obtained:

$$\frac{\partial l(\theta)}{\partial \mu} = -\sum_{i=1}^{n}-\frac{(y_i-\mu)}{\sigma^2} + 2\sum_{i=1}^{n}\frac{2\alpha b_i}{\sigma(1+b_i^2)}$$

$$\frac{\partial l(\theta)}{\partial \sigma} = -\frac{n}{\sigma} - \sum_{i=1}^{n}-\frac{(y_i-\mu)^2}{\sigma^3} + 2\sum_{i=1}^{n}\frac{2\alpha(y_i-\mu)b_i}{\sigma^2(1+b_i^2)}$$

$$\frac{\partial l(\theta)}{\partial \alpha} = -\frac{n(16\alpha+12\alpha^3)}{4+8\alpha^2+3\alpha^4} + 2\sum_{i=1}^{n}-\frac{2(y_i-\mu)b_i}{\sigma(1+b_i^2)}$$

where, $b_i = \left(1 - \frac{\alpha(y_i-\mu)}{\sigma}\right)$

The solutions of the above system of likelihood equations gives the maximum likelihood estimator for $\theta = (\alpha, \mu, \sigma)$ which can be obtained by numerically maximizing eqn. (26) with respect to the parameters $\theta = (\alpha, \mu, \sigma)$. The derivation of the observed information matrix is also obtained using numerical procedures. Initial values for those procedures can be obtained although the moment estimators.

The log-likelihood function based on a single observation $Y$ for the parameters $\theta = (\alpha, \mu, \sigma)$, is given in Appendix F.

The variance-covariance matrix of the MLEs can be obtained by taking the inverse of the Fisher information matrix (**I**) as given in Appendix F.

## 6. Real life applications: comparative data fitting

Here we have considered two datasets which are related N latitude degrees in 69 samples from world lakes, which appear in Column 5 of the Diversity data set in website:

http://users.stat.umn.edu/sandy/cours es/8061/datasets/lakes.lsp; and the body mass index (BMI) of 202 Australian athletes (Cook and Weisberg, 1994).

We then compared the proposed distribution $BASN_2(\alpha,\mu,\sigma)$ with the normal distribution $N(\mu,\sigma^2)$, the logistic distribution $LG(\mu,\beta)$, the Laplace distribution $La(\mu,\beta)$, the skew-normal distribution $SN(\lambda,\mu,\sigma)$ of Azzalini (1985), the skew-logistic distribution $SLG(\lambda,\mu,\beta)$ of Wahed and Ali (2001), the skew-Laplace distribution $SLa(\lambda,\mu,\beta)$ of Nekoukhou and Alamatsaz (2012), the alpha-skew-normal distribution $ASN(\alpha,\mu,\sigma)$ of Elal-Olivero (2010), the alpha-skew-Laplace distribution $ASLa(\alpha,\mu,\beta)$ of Harandi and Alamatsaz (2013), the alpha-skew-logistic distribution $ASLG(\alpha,\mu,\beta)$ of Hazarika and Chakraborty (2014), the alpha-beta-skew-normal distribution $ABSN(\alpha,\beta,\mu,\sigma)$ and beta-skew-normal distribution $BSN(\beta,\mu,\sigma)$ of Shafiei et al. (2016).

Using GenSA package in R, the MLE of the parameters are obtained by using numerical optimization routine. AIC and BIC are used for model comparison.
Table 1 and Table 2 shows the MLE's, log-likelihood, AIC and BIC of the above mentioned distributions. The graphical representation of the results taking only the top three competitors for the proposed model are given in figure 10 and figure 11.

**Table 1:** MLE's, log-likelihood, AIC and BIC for N latitude degrees in 69 samples from world lakes.

| Parameters → <br> Distribution ↓ | $\mu$ | $\sigma$ | $\lambda$ | $\alpha$ | $\beta$ | $\log L$ | AIC | BIC |
|---|---|---|---|---|---|---|---|---|
| $N(\mu,\sigma^2)$ | 45.165 | 9.549 | -- | -- | -- | -253.599 | 511.198 | 515.666 |
| $LG(\mu,\beta)$ | 43.639 | -- | -- | -- | 4.493 | -246.645 | 497.290 | 501.758 |
| $SN(\lambda,\mu,\sigma)$ | 35.344 | 13.70 | 3.687 | -- | -- | -243.036 | 492.072 | 498.774 |
| $BSN(\beta,\mu,\sigma)$ | 54.47 | 5.52 | -- | -- | 0.74 | -242.530 | 491.060 | 497.760 |
| $SLG(\lambda,\mu,\beta)$ | 36.787 | -- | 2.8284 | -- | 6.417 | -239.053 | 490.808 | 490.808 |
| $La(\mu,\beta)$ | 43.00 | -- | -- | -- | 5.895 | -239.248 | 482.496 | 486.964 |
| $ASLG(\alpha,\mu,\beta)$ | 49.087 | -- | -- | 0.861 | 3.449 | -237.351 | 480.702 | 487.404 |
| $SLa(\lambda,\mu,\beta)$ | 42.30 | -- | 0.255 | -- | 5.943 | -236.900 | 479.799 | 486.501 |
| $ASLa(\alpha,\mu,\beta)$ | 42.3 | -- | -- | -0.220 | 5.440 | -236.079 | 478.159 | 484.861 |
| $ASN(\alpha,\mu,\sigma)$ | 52.147 | 7.714 | -- | 2.042 | -- | -235.370 | 476.739 | 483.441 |
| $ABSN(\alpha,\beta,\mu,\sigma)$ | 47.69 | 7.15 | -- | 1.72 | -0.37 | -230.770 | 469.530 | 478.480 |
| $BASN_2(\alpha,\mu,\sigma)$ | 54.265 | 6.559 | -- | 1.994 | -- | **-226.228** | **458.455** | **465.158** |

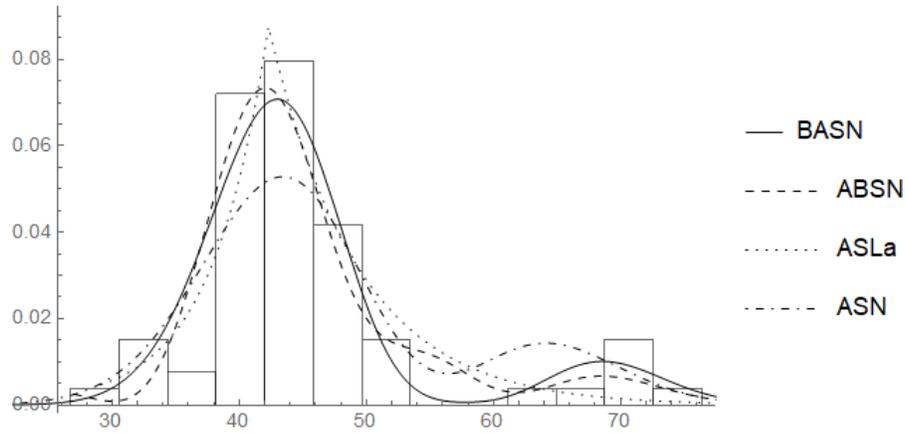

**Figure 10:** Plots of observed and expected densities of some distributions for N latitude degrees in 69 samples from world lakes.

**Table 2:** MLE's, log-likelihood, AIC and BIC for body mass index (BMI) of 202 Australian athletes.

| Parameters → Distribution ↓ | $\mu$ | $\sigma$ | $\lambda$ | $\alpha$ | $\beta$ | $\log L$ | AIC | BIC |
|---|---|---|---|---|---|---|---|---|
| $N(\mu,\sigma^2)$ | 22.956 | 2.857 | -- | -- | -- | -498.668 | 1001.336 | 1007.953 |
| $La(\mu,\beta)$ | 22.749 | -- | -- | -- | 2.123 | -494.08 | 992.16 | 998.7765 |
| $BSN(\beta,\mu,\sigma)$ | 22.528 | 2.694 | -- | -- | -0.058 | -492.88 | 991.76 | 1001.685 |
| $ASLa(\alpha,\mu,\beta)$ | 22.350 | -- | -- | -0.14 | 2.07 | -492.601 | 991.202 | 1001.127 |
| $SLa(\lambda,\mu,\beta)$ | 22.350 | -- | 0.865 | -- | 2.084 | -492.461 | 990.922 | 1000.847 |
| $LG(\mu,\beta)$ | 22.787 | -- | -- | -- | 1.529 | -491.462 | 986.924 | 993.5405 |
| $SN(\lambda,\mu,\sigma)$ | 19.969 | 4.133 | 2.313 | -- | -- | -490.099 | 986.198 | 996.1228 |
| $ASLG(\alpha,\mu,\beta)$ | 21.933 | -- | -- | -0.201 | 1.475 | -489.094 | 984.188 | 994.1128 |
| $ASN(\alpha,\mu,\sigma)$ | 24.834 | 2.653 | -- | 0.994 | -- | -488.69 | 983.38 | 993.3048 |
| $ABSN(\alpha,\beta,\mu,\sigma)$ | 23.998 | 2.853 | -- | 0.817 | -0.131 | -486.743 | 981.486 | 994.7191 |
| $SLG(\lambda,\mu,\beta)$ | 20.717 | -- | 1.401 | -- | 1.975 | -487.311 | 980.622 | 990.5468 |
| $BASN_2(\alpha,\mu,\sigma)$ | 26.482 | 2.706 | -- | 0.971 | -- | **-484.773** | **975.546** | **985.4708** |

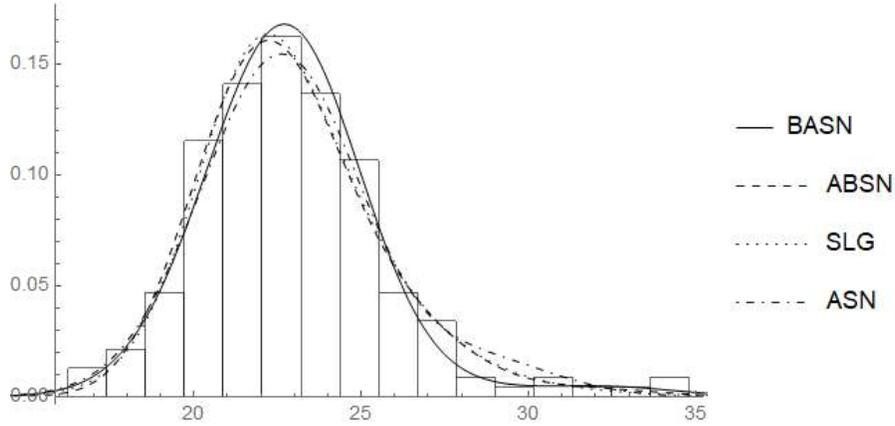

**Figure 11:** Plots of observed and expected densities of some distributions for body mass index (BMI) of 202 Australian athletes.

It is found from the Table 1 and 2 that the proposed Balakrishnan-alpha-skew-normal distribution $BASN_2(\alpha,\mu,\sigma)$ provides best fit to the data set in terms of AIC and BIC. The plots of observed and expected densities presented in Figure 10 and 11, also confirm our findings.

**Remark 11:** The observed variance-covariance matrix of the MLEs of the parameters $\theta = (\alpha,\mu,\sigma)$ of $BASN_2(\alpha,\mu,\sigma)$ distribution in example 1 and example 2 are obtained respectively, as

$$Var-Cov(\hat{\theta}) = \begin{pmatrix} 0.6995 & 0.1615 & 0.1395 \\ 0.1615 & 0.1350 & -0.00455 \\ 0.1395 & -0.00455 & 0.1349 \end{pmatrix} \& \ Var-Cov(\hat{\theta}) = \begin{pmatrix} 0.1371 & 0.0327 & 0.0337 \\ 0.0327 & 0.0199 & 0.00192 \\ 0.0337 & 0.00192 & 0.0175 \end{pmatrix}$$

### 6.1. Likelihood Ratio Test

Furthermore, since $N(\mu,\sigma^2)$ and $BASN_2(\alpha,\mu,\sigma)$ distributions are nested models, the likelihood ratio (LR) test is used to differentiate between them. The LR test is carried out to test the following null hypothesis $H_0 : \alpha = 0$, that is the sample is drawn from $N(\mu,\sigma^2)$; against the alternative $H_1 : \alpha \neq 0$, that is the sample is drawn from $BASN_2(\alpha,\mu,\sigma)$.

The values of LR test statistic for the data set I and II are respectively 54.742 and 27.79. Both of which exceed the 99% critical value, that is, 6.635. Thus there is evidence in support of the alternative hypothesis that is, the sampled data comes from $BASN_2(\alpha,\mu,\sigma)$, not from $N(\mu,\sigma^2)$.

## 7. Concluding remarks

In this study a new alpha-skew-normal distribution with one parameter which has both unimodal as well as bimodal shapes is constructed and some of its properties are studied. The bathtub shaped failure rate function is seen in the half $BASN_2(\alpha, \mu, \sigma)$ distribution. Some extensions of the proposed distribution with some of their particulars cases are presented. Our findings adequately supported the proposed $BASLG_2(\alpha, \mu, \beta)$ distribution as the best fitted one to the datasets under consideration in terms of AIC and BIC. The plots of observed and expected densities presented also confirm our findings.

## Appendix

### A: Proof of Proposition 1

$$F_Z(z) = = \frac{1}{C_2(\alpha)(2+\alpha^2)^2} \int_{-\infty}^{\infty} \left\{\alpha^4 z^4 - 4\alpha^3 z^3 + 8\alpha^2 z^2 - 8\alpha z + 4\right\} \varphi(z)\, dz$$

$$= \frac{\alpha^4[-\{z(3+z^2)\varphi(z)\} + 3\Phi(z)] - 4\alpha^3[-\{(2+z^2)\varphi(z)\}] + 8\alpha^2\{-z\varphi(z)+\Phi(z)\} - 8\alpha(-\varphi(z)) + 4\Phi(z)}{C_2(\alpha)(2+\alpha^2)^2}$$

$$= \frac{\varphi(z)[-3\alpha^4 z - \alpha^4 z^3 + 8\alpha^3 + 4\alpha^3 z^2 - 8\alpha^2 z + 8\alpha]}{C_2(\alpha)(2+\alpha^2)^2} + \frac{\Phi(z)[3\alpha^4 + 8\alpha^2 + 4]}{(3\alpha^4 + 8\alpha^2 + 4)}$$

$$= \Phi(z) + \frac{\alpha[8 - 8\alpha z + 4\alpha^2(2+z^2) - \alpha^3 z(3+z^2)]}{C_2(\alpha)(2+\alpha^2)^2} \varphi(z)$$

### B: Proof of Proposition 2

When *n* is even; $\quad E(Z^n) = \int_{-\infty}^{\infty} \frac{z^n}{C_2(\alpha)} \left\{ \frac{(1-\alpha z)^2 + 1}{2+\alpha^2} \right\}^2 \varphi(z)\, dz$

$$= \frac{1}{C_2(\alpha)(2+\alpha^2)^2} \int_{-\infty}^{\infty} \left( \alpha^4 z^{n+4} - 4\alpha^3 z^{n+3} + 8\alpha^2 z^{n+2} - 8\alpha z^{n+1} + 4z^n \right) \varphi(z) dz$$

$$= \frac{1}{C_2(\alpha)(2+\alpha^2)^2} \left[ \alpha^4 \int_{-\infty}^{\infty} z^{n+4} \varphi(z) dz + 8\alpha^2 \int_{-\infty}^{\infty} z^{n+2} \varphi(z) dz + 4 \int_{-\infty}^{\infty} z^n \varphi(z) dz \right]$$

using the result of $n^{th}$ order moment of normal distribution in the above equation we get the result in (9). Similarly, when $n$ is odd the same can be obtained with the help of following equation.

$$E(Z^n) = \frac{1}{C_2(\alpha)(2+\alpha^2)^2} \left[ -4\alpha^3 \int_{-\infty}^{\infty} z^{n+3} \varphi(z) dz - 8\alpha \int_{-\infty}^{\infty} z^{n+1} \varphi(z) dz \right]$$

**C: Proof of Proposition 3**

$$M_Z(t) = \frac{1}{C_2(\alpha)(2+\alpha^2)^2} \int_{-\infty}^{\infty} e^{tz} \left[ \alpha^4 z^4 - 4\alpha^3 z^3 + 8\alpha^2 z^2 - 8\alpha z + 4 \right] \varphi(z) dz \tag{A1}$$

Again, $\int_{-\infty}^{\infty} \varphi(z) e^{tz} dz = e^{\frac{t^2}{2}} = M_X(t)$, i.e., mgf of standard normal variable

$$\int_{-\infty}^{\infty} z \varphi(z) e^{tz} dz = t M_X(t), \quad \int_{-\infty}^{\infty} z^2 \varphi(z) e^{tz} dz = t^2 M_X(t) + M_X(t)$$

$$\int_{-\infty}^{\infty} z^3 \varphi(z) e^{tz} dz = t^3 M_X(t) + 3t M_X(t) \text{ and } \int_{-\infty}^{\infty} z^4 \varphi(z) e^{tz} dz = t^4 M_X(t) + 6t^2 M_X(t) + 3M_X(t)$$

Now, applying the above results in (A1) we get the expression in (10).

**D: Proof of Proposition 4**

$$F(z) = \frac{1}{(2+\alpha^2)(2+3\alpha^2)} \int_{-\infty}^{z} \left[ \alpha^4 z^4 + 8\alpha^2 z^2 + 4 \right] \left[ \frac{1}{\sqrt{2\pi}} e^{-\frac{z^2}{2}} \right] dz$$

$$= \frac{\alpha^4 [-z(3+z^2)\varphi(z) + 3\Phi(z)] + 8\alpha^2 [-z\varphi(x) + \Phi(z)] + 4\Phi(z)}{C_2(\alpha)(2+\alpha^2)^2}$$

On simplifying the above equation we get the result in (12).

**E: Proof of Proposition 4**

$$M_Z(t) = \frac{1}{C_2(\alpha)(2+\alpha^2)^2} \int_{-\infty}^{\infty} e^{tz} \left[ \alpha^4 z^4 + 8\alpha^2 z^2 + 4 \right] \varphi(z) dz \tag{A2}$$

Again, $\int_{-\infty}^{\infty} \varphi(z) e^{tz} dz = e^{\frac{t^2}{2}} = M_X(t)$, i.e., mgf of standard normal variable

$$\int_{-\infty}^{\infty} z^2 \varphi(z) e^{tz} dz = t^2 M_X(t) + M_X(t) \text{ and } \int_{-\infty}^{\infty} z^4 \varphi(z) e^{tz} dz = t^4 M_X(t) + 6t^2 M_X(t) + 3M_X(t)$$

Now, applying the above results in (A2), we get the expression in (13).

# F: Log-Likelihood and Fisher Information Matrix

$$l(\theta;y) = 2\log\left[\left\{1-\alpha\left(\frac{y-\mu}{\sigma}\right)\right\}^2 + 1\right] - \log(2+\alpha^2) - \log\sigma - \log(2+3\alpha^2) - \frac{1}{2}\log(2\pi) - \frac{1}{2}\left(\frac{y-\mu}{\sigma}\right)^2$$

*Score functions:*

The first-order partial derivatives of $l(\theta;y)$ are:

$$\frac{\partial l(\theta;y)}{\partial \mu} = \frac{(y-\mu)}{\sigma^2} + \frac{4\alpha b}{\sigma(1+b^2)}$$

$$\frac{\partial l(\theta;y)}{\partial \sigma} = -\frac{n}{\sigma} + \frac{(y-\mu)^2}{\sigma^3} + \frac{4\alpha(y-\mu)b}{\sigma^2(1+b^2)}$$

$$\frac{\partial l(\theta;y)}{\partial \alpha} = -\frac{n(16\alpha + 12\alpha^3)}{4+8\alpha^2+3\alpha^4} - \frac{4(y-\mu)b}{\sigma(1+b^2)}$$

The second-order partial derivatives of $l(\theta;y)$ are:

$$\frac{\partial^2 l(\theta;y)}{\partial \mu^2} = -\frac{1}{\sigma^2} + 2\left(\frac{2\alpha^2}{\sigma^2(1+b^2)} - \frac{4\alpha^2 b^2}{\sigma^2(1+b^2)^2}\right)$$

$$\frac{\partial^2 l(\theta;y)}{\partial \sigma^2} = \frac{n}{\sigma^2} - \frac{3(y-\mu)^2}{\sigma^4} + 2\left(\frac{2\alpha^2(y-\mu)^2}{\sigma^4(1+b^2)} - \frac{4\alpha^2(y-\mu)^2 b^2}{\sigma^4(1+b^2)^2} - \frac{4\alpha(y-\mu)b}{\sigma^3(1+b^2)}\right)$$

$$\frac{\partial^2 l(\theta;y)}{\partial \alpha^2} = -n\left(-\frac{(16\alpha+12\alpha^3)^2}{(4+8\alpha^2+3\alpha^4)^2} + \frac{16+36\alpha^2}{4+8\alpha^2+3\alpha^4}\right) + 2\left(\frac{2(y-\mu)^2}{\sigma^2(1+b^2)} - \frac{4(y-\mu)^2 b^2}{\sigma^2(1+b^2)^2}\right)$$

$$\frac{\partial^2 l(\theta;y)}{\partial \mu \partial \sigma} = -\frac{2(y-\mu)}{\sigma^3} + \frac{4\alpha^2(y-\mu)}{\sigma^3(1+b^2)} - \frac{8\alpha^2(y-\mu)b^2}{\sigma^3(1+b^2)^2} - \frac{4\alpha b}{\sigma^2(1+b^2)}$$

$$\frac{\partial^2 l(\theta;y)}{\partial \mu \partial \alpha} = -\frac{4\alpha(y-\mu)}{\sigma^2(1+b^2)} + \frac{8\alpha(y-\mu)b^2}{\sigma^2(1+b^2)^2} + \frac{4b}{\sigma(1+b^2)}$$

$$\frac{\partial^2 l(\theta;y)}{\partial \sigma \partial \alpha} = -\frac{4\alpha(y-\mu)^2}{\sigma^3(1+b^2)} + \frac{8\alpha(y-\mu)^2 b^2}{\sigma^3(1+b^2)^2} + \frac{4(y-\mu)b}{\sigma^2(1+b^2)}$$

where, $b = \left(1 - \frac{\alpha(y-\mu)}{\sigma}\right)$

$$I = \begin{bmatrix} E\left(-\frac{\partial^2 l(\theta;y)}{\partial \mu^2}\right) & E\left(-\frac{\partial^2 l(\theta;y)}{\partial \mu \partial \sigma}\right) & E\left(-\frac{\partial^2 l(\theta;y)}{\partial \mu \partial \alpha}\right) \\ E\left(-\frac{\partial^2 l(\theta;y)}{\partial \sigma \partial \mu}\right) & E\left(-\frac{\partial^2 l(\theta;y)}{\partial \sigma^2}\right) & E\left(-\frac{\partial^2 l(\theta;y)}{\partial \sigma \partial \alpha}\right) \\ E\left(-\frac{\partial^2 l(\theta;y)}{\partial \alpha \partial \mu}\right) & E\left(-\frac{\partial^2 l(\theta;y)}{\partial \alpha \partial \sigma}\right) & E\left(-\frac{\partial^2 l(\theta;y)}{\partial \alpha^2}\right) \end{bmatrix}^{-1}$$

Where, $E\left(-\frac{\partial^2 l(\theta;y)}{\partial \mu^2}\right) \approx -\frac{\partial^2 l(\theta;y)}{\partial \mu^2}\bigg|_{\mu=\hat{\mu}}$, $E\left(-\frac{\partial^2 l(\theta;y)}{\partial \mu \partial \alpha}\right) \approx -\frac{\partial^2 l(\theta;y)}{\partial \mu \partial \alpha}\bigg|_{\mu=\hat{\mu}, \alpha=\hat{\alpha}}$ etc.